\documentclass[a4paper,11pt]{article}
\usepackage{lmodern}
\usepackage{amsmath,amsfonts,amssymb,mathrsfs,charter,amsthm}

\title{Orthogonality to matrix subspaces, and a distance formula}
\author{Priyanka Grover\\
  \multicolumn{1}{p{.7\textwidth}}{\centering\emph{Theoretical Statistics and Mathematics Unit, Indian Statistical Institute, Delhi Centre, 7, S.J.S. Sansanwal Marg, New Delhi-110016, India\\ Email:  pgrover8r@isid.ac.in}}}
\date{}

\newcommand{\mat}{\mathbb{M}_{n}(\mathbb C)}

\newcommand{\W}{\mathscr W}
\newcommand{\C}{\mathbb{C}}
\newcommand{\tr}{{\operatorname{tr\ }}}
\newcommand{\D}{\mathbb D}

\newcommand{\R}{\mathbb{R}}
\newcommand{\h}{{\operatorname{Re}}}

\newcommand{\dist}{{\operatorname{dist}}}
\newcommand{\diag}{{\operatorname{diag}}}
\newcommand{\conv}{\mathop{{\rm conv}}}
\newtheorem{thm}{Theorem}

\newtheorem{prop}{Proposition}

\newcommand{\A}{\mathscr{B}}

\begin{document}

{\setlength{\baselineskip}%
{1\baselineskip}
\maketitle
\begin{abstract}
We obtain a necessary and sufficient condition for a matrix $A$ to be Birkhoff-James orthogonal to any subspace $\mathscr W$ of $\mat$. Using this we obtain an expression for the distance of $A$ from any unital $C^*$ subalgebra of $\mat$.
\end{abstract}


}

\section{Introduction}
Let $\mat$ be the space of $n \times n$ complex matrices and let $\W$ be any subspace of $\mat$.  For any $A\in \mat$, let $$\|A\|=\max_{ x\in \C^n, \|x\|=1} \|Ax\|$$ be the \emph{operator norm} of $A$. Then $A$ is said to be \emph{(Birkhoff-James) orthogonal} to $\W$ if \begin{equation}\|A+W\|\geq \|A\|  \text{ for all } W \in \W.\end{equation}
The space $\mat$ is a complex Hilbert space under the inner product $\langle A, B\rangle_c = \tr (A^* B)$ and a real Hilbert space under the inner product $\langle A, B\rangle_r =\h\ \tr (A^* B)$. Let $\W^{\perp}$ be the orthogonal complement of $\W$, where the orthogonal complement is with respect to the usual Hilbert space orthogonality in $\mat$ with the inner product $\langle \cdot, \cdot \rangle_c$ or $\langle \cdot, \cdot \rangle_r$, depending upon whether $\W$ is a real or complex subspace.  Note that if $A\in \W^{\perp}$ such that $\tr(A^* A)=\|A\|^2 $, then $A$ is orthogonal to $\W$.

Bhatia and $\check{\text S}$emrl \cite{semrl} obtained an interesting characterisation of orthogonality when $\W=\C B$, where $B$ is any matrix in $\mat$. They showed that $A$ is orthogonal to $\C B$ if and only if there exists a unit vector $x$ such that $\|Ax\|=\|A\|$
and $\langle Ax,Bx \rangle=0$. In other words, $A$ is orthogonal to $\C B$ if and only if there exists a positive semidefinite matrix $P$ of rank one such that $\tr P=1, \ \tr A^* AP=\|A\|^2$ and $AP\in (\C B)^{\perp}.$ Such positive semidefinite matrices with trace 1 are called \emph{density matrices}. We use the notation $P\geq 0$ to mean $P$ is positive semidefinite.

Let $\W=\D_n(\R),$ the subspace of all diagonal matrices with real entries, and let $A$ be any Hermitian matrix. Then $A$ is called \emph{minimal} if $\|A+D\|\geq \|A\|$ for all $D\in \D_n(\R)$. Andruchow, Larotonda, Recht, and Varela \cite[Theorem 1]{andruchow} showed that a Hermitian matrix $A$ is minimal if and only if there exists a density matrix $P$ such that  $\ PA^2=\|A\|^2 P$ and all diagonal entries of $PA$ are zero. In our notation, $A$ is minimal is same as saying that $A$ is orthogonal to the subspace $\D_n(\R)$. If $A$ is Hermitian, then note that $A$ is orthogonal to $\D_n(\R)$ if and only if $A$ is orthogonal to $\D_n(\C)$. Now $\D_n(\C)^{\perp}$ is the subspace of all matries such that their diagonal entries are zero. The condition  $\ PA^2=\|A\|^2 P$ is same as $A^2 P=\|A\|^2 P$ and diagonal entries of $PA$ are same as diagonal entries of $AP$. Therefore  Theorem 1 in \cite{andruchow} can be interpreted as follows. A Hermitian matrix $A$ is orthogonal to $\D_n(\C)$ if any only if $ A^2 P=\
|A\|^2 P$ and  $AP\in \D_n(\C)^{\perp}$. The following theorem is a generalization of  this result as well as Bhatia-$\check{\text S}$emrl theorem.

\begin{thm}\label{thm1}
Let $A\in \mat$ and let $m(A)$ be the multiplicity of the maximum singular value $\|A\|$ of $A$. Let $\W$ be any (real or complex) subspace of $\mat.$ Then $A$ is orthogonal to $\W$ if and only if there exists a density matrix $P$ of complex rank at most $m(A)$ such that $A^*AP=\|A\|^2 P$ and $AP \in \W^{\perp}$.
(If rank $P=\ell$, then $P$ has the form $P=\displaystyle\sum_{i=1}^{\ell} t_i v_{(i)}v_{(i)}^*$ 
where $v_{(i)}$ are unit vectors such that $A^*A v_{(i)}=\|A\|^2 v_{(i)}$ and $t_i$ are such that $0\leq t_i\leq 1$ and $\displaystyle\sum_{i=1}^{\ell} t_i=1$.)
\end{thm}

Here, $m(A)$ is the  best possible upper bound on rank $P$. This has been illustrated later in Remark 4 in Section 4. 
When $\W=\C B$, the above theorem says that  $A$ is orthogonal to $\C B$ if and only if there exists a $P\geq 0$ of the form $P=\displaystyle{\sum_{i=1}^{\ell} }t_i v_{(i)}v_{(i)}^*$ such that $\|v_{(i)}\|=1$, $A^*Av_{(i)}=\|A\|^2 v_{(i)}$ and $\displaystyle\sum_{i=1}^{\ell}t_i\langle B^* A v_{(i)}, v_{(i)}\rangle=0$. By the Hausdorff-Toeplitz theorem, we get a unit vector $v$ such that $A^*A v=\|A\|^2v$ and $\langle B^*  A v, v\rangle=0$. The first condition is stronger than that in \cite[Theorem 1.1]{semrl}. 

Let $\dist(A, \W)$ denote the distance of a matrix $A$ from the subspace $\W$, defined as
$$\dist(A,\W)=\min \left\{\|A-W\|: W\in \W\right\}.$$
Audenaert \cite{audenaert} showed that when $\W=\C I$, then
\begin{equation}
\dist(A,\C I)^2=\max{\left\{\tr(A^*AP)-|\tr(AP)|^2: P\geq 0,\tr P=1\right\}}.\label{dist to CI}
\end{equation}
Further the maximisation over $P$ on the right hand side of \eqref{dist to CI} can be restricted to density matrices of rank 1. The quantity $\tr(A^*AP)-|\tr(AP)|^2$ is called the \emph{variance} of $A$ with respect to the density matrix $P$. 
Bhatia and Sharma \cite{sharma} showed that if $\Phi: \mat \rightarrow \mathbb M_k(\C)$ is any positive unital linear map, then 
$$\Phi(A^*A)-\Phi(A)^* \Phi(A)\leq \dist(A,\C I)^2.$$ By choosing $\Phi(A)=\tr(AP)$ for different density matrices $P$, they obtained various interesting bounds on $\dist(A,\C I)^2.$

 It would be interesting to have a generalisation of \eqref{dist to CI} with  $\C I$ replaced by any unital $C^*$ subalgebra of $\mat$. (This problem has also been raised by M. Rieffel in \cite{rieffel}.) Let $\A$ be any unital $C^*$ subalgebra of $\mat$. Let $\mathcal C_{\A}: \mat \rightarrow \A$ denote the projection of $\mat$ onto $\A$. We note that $\mathcal C_{\A}$ is a bimodule map:
\begin{equation}
\mathcal C_{\A}(BX)=B \mathcal C_{\A}(X) \text{ and } \mathcal C_{\A}(XB)=\mathcal C_{\A}(X) B \text{ for all } B\in \A, X\in \mat.\label{bimodule}
\end{equation}
In particular, when $\A$ is the subalgebra of block diagonal matrices, 
the matrix $\mathcal C_{\A}(X)$ is called a \emph{pinching} of $X$ and is denoted by $\mathcal C(X)$.
It is defined as follows. 
If $X=\left[\begin{array}{cccc} X_{11} &\cdots & X_{1k}\\ X_{21} & \cdots & X_{2k}\\ 
\vdots & \vdots & \vdots\\ X_{k1}&\cdots & X_{kk}\end{array}\right]$ then 
\begin{equation}
\mathcal C(X)=
\left[\begin{array}{cccc} X_{11} &\ & \ &\ \\ \ & X_{22} & \ &\ \\  \ & \ &\ddots & \ \\ 
\ &\ &\ & X_{kk}\end{array}\right].\label{pinching}
\end{equation} 
Properties of pinchings are studied in detail in \cite{bhatiama} and \cite{bhatiapdm}.

 Our next result provides a generalisation of \eqref{dist to CI} for distance of $A$ to any unital $C^*$ subalgebra of $\mat$. 
\begin{thm}\label{thm2}
Let $\A$ be any unital $C^*$ subalgebra of $\mat$. Let $\mathcal C_{\A}: \mat \rightarrow \A$ denote the projection of $\mat$ onto $\A$.

Then 
\begin{equation}
\dist(A, \A)^2=\max\left\{\tr\left(A^*AP-\mathcal C_{\A}(AP)^*\ \mathcal C_{\A}(AP)\ \mathcal C_{\A}(P)^{-1}\right): P\geq 0,\tr P=1 \right\},\label{dist to any subalgebra}
\end{equation}
where $\mathcal C_{\A}(P)^{-1}$ denotes the Moore-Penrose inverse of $\mathcal C_{\A}(P)$. The maximum on the right hand side of \eqref{dist to any subalgebra} can be restricted to rank $P\leq m(A)$.
\end{thm}

We prove Theorem 1 using ideas of subdifferential calculus. A brief summary of these is given in Section 2. The proofs are given in Section 3.

\section{Preliminaries}

Let $X$ be a complex Hilbert space. Let $f:X\rightarrow \R$ be a convex function. 
Then the \emph{subdifferential} of $f$ at any point $x\in X$, denoted by $\partial f(x)$, 
is the set of $v^*\in X^*$ such that
\begin{equation}f(y)-f(x)\geq \h\ v^*(y-x) \text{ for all } y\in X.\label{subdifferential}
\end{equation} 
It follows from \eqref{subdifferential} that $f$ is minimized at $x$ if and only if 
$0\in \partial f(x)$. 

We use an idea similar to the one in \cite[Theorem 2.1]{tbgrover}. 
Let $f(W)=\|A+W\|$. This is the composition of two functions namely $W\rightarrow A+W$ from $\W$ 
into $\mat$ and $T\rightarrow \|T\|$ from $\mat$ into $\R_+$. 
Thus we need to find subdifferentials of composition maps. For that we need a chain rule.

\begin{prop}\label{p3}
Let $X,Y$ be any two Hilbert spaces. Let $g:Y\rightarrow \R$ be a convex function. 
Let $S:X\rightarrow Y$ be a linear map and let $L:X\rightarrow Y$ be the affine map 
defined by $L(x)=S(x)+y_0$, for some $y_0\in Y$. Then
\begin{equation}
\partial(g\circ L)(x)=S^* \partial g(L(x)),\label{composition}
\end{equation}
where $S^*$ is the adjoint of $S$ defined as 
$$\langle S^*(y),x\rangle=\langle y, S(x)\rangle \text{ for all } x\in X \text{ and } y\in Y.$$ 
\end{prop}

In our setting, $g$ is the map $T\rightarrow \|T\|$. The subdifferential of this map has been calculated by Watson \cite{watson}. 

\begin{prop}\label{p2}
 Let $A\in \mat$. Then
\begin{equation}
\partial \|A\|=\conv\{uv^*:  \|u\|=\|v\|=1, Av=\|A\|u \},\label{partial}
\end{equation}
where $\conv D$ denotes the  convex hull of a set $D$.
\end{prop}

These elementary facts can be found in \cite{hiriart}. In this book the author deals with convex functions $f: \R^n\rightarrow \R$. The same proofs can be extended to functions $f:X\rightarrow \R$, where $X$ is any Hilbert space.

\section{Proofs}
\textit{Proof of Theorem \ref{thm1}} \quad
Suppose there exists a positive semidefinite $P$ with $\tr P=1$ such that 
$A^*AP=\|A\|^2 P$ and $AP \in \W^{\perp}$. Then for any $W\in\W$
\begin{eqnarray*}
\|A+W\|^2&=&\|(A+W)^*(A+W)\| \nonumber\\
&=& \|A^*A + W^* A + A^* W +W^* W\|. \nonumber
\end{eqnarray*}
Now for any $T\in \mat$,
\begin{equation}
\|T\|=\sup_{\|X\|_1=1} |\tr(T  X)|,\label{supnorm}
\end{equation} where $\|\cdot\|_1$ denotes the trace norm. So,
\begin{eqnarray}
\|A+W\|^2 &\geq& |\tr (A^*A P+ W^* A P+ A^* W P+W^* W P)| \nonumber \\
&\geq& \h \ \tr (A^*A P+ W^* A P+ A^* W P+W^* W P) .\label{trace inequality}
\end{eqnarray}
Since $AP\in \W^{\perp}$, we have $\h\ \tr(A^* W P)=\h\ \tr (W^*AP)=0$. 
 The matrices $W^* W$ and $P$ are positive semidefinite, therefore $\tr (W^*WP)\geq 0$ and by our assumption, $\tr(A^*AP)=\|A\|^2$. Using these in \eqref{trace inequality} we get that $\|A+W\|^2\geq \|A\|^2$.

Conversely, suppose \begin{equation}
\|A+W\|\geq \|A\|\text{ for all } W\in \W.\label{mainresult}
\end{equation}
Let $S: \W\rightarrow \mat$ be the inclusion map. Then $S^*:\mat \rightarrow \W$ is the 
projection onto the subspace $\W$.  Let $L:\W \rightarrow \mat$ be the map defined as
$$L(W)=A+S(W).$$
Let $g: \mat \rightarrow \R$ be the map taking an $n\times n$ matrix $W$ to $\|W\|$.
Then \eqref{mainresult} can be rewritten as 
$$(g\circ L)(W)\geq (g\circ L)(0),$$
that is, $g\circ L$ is minimized at $0$. Therefore $0 \in \partial (g\circ L)(0)$. 
Using Proposition \ref{p3}, we get
\begin{equation}
0\in S^*\partial \|A\|.\label{zeroinsubdifferential}
\end{equation}
By Proposition \ref{p2},
\begin{equation}
S^*\partial \|A\|=\rm{conv}\left\{S^*(uv^*): \|u\|=\|v\|=1, Av=\|A\|u\right\}.\label{theset}
\end{equation}

From \eqref{zeroinsubdifferential} and \eqref{theset} it follows that there exist unit vectors $u_{(i)}, v_{(i)}$ such that $Av_{(i)}=\|A\|u_{(i)}$ and numbers $t_i$ such that $0\leq t_i\leq 1$, $\sum t_i=1$ and 
\begin{equation}
S^*\left(\sum t_i u_{(i)} v_{(i)}^*\right)=0.\label{condition}
\end{equation}
Let $P=\sum t_i v_{(i)} v_{(i)}^*$. Then $P\geq 0$ and $\tr P=1$. Note that 
\begin{eqnarray*}
AP&=&\sum t_i Av_{(i)} v_{(i)}^*\\
&=&\|A\|\sum t_i u_{(i)} v_{(i)}^*.
\end{eqnarray*}
So, from \eqref{condition} we get $S^*(AP)=0$, that is, $AP\in \W^{\perp}$.
Since each $v_{(i)}$ is a right singular vector for $A$, we have $A^*A v_{(i)}=\|A\|^2  v_{(i)}$. Using this we obtain
\begin{eqnarray}
A^*AP&=&\sum t_i A^*A v_{(i)} v_{(i)}^*\nonumber\\
&=& \sum t_i \|A\|^2 v_{(i)} v_{(i)}^*\nonumber\\
&=&\|A\|^2 P. \label{consequence1}
\end{eqnarray}

Now let $m(A)=k$. We now show that if $P$ satisfies \eqref{consequence1}, then rank $P\leq k$. First note that $A^*A$ and $P$ commute and therefore can be diagonalised simultaneously. So we can assume $A^*A$ and $P$ in \eqref{consequence1} to be diagonal matrices. By hypothesis $k$ of the diagonal entries of $A^*A$ are equal to $\|A\|^2$. Let $A^*A=\left[
\begin{array}{ccccccc}
\|A\|^2 & \         &\  &\      &\            &\         &\      \\
\            & \ddots&\ &  \     &\            &\         &\      \\
\            & \          &\|A\|^2&\            &\         &\      \\
\            &\           &\           &s_{k+1}^2&\         &\      \\
\            &\           &\           &\            &\ddots&\      \\
\            &\           &\           &\            &\         &s_n^2
\end{array}
\right],$ where $s_j<\|A\|$ for all $k+1\leq j\leq n$. 
If $P=\left[\begin{array}{ccc} p_1 & \ &\ \\ \ &\ddots&\ \\ \ &\ &p_n \end{array}\right],$ 
then from \eqref{consequence1} we obtain 
$$(s_j^2-\|A\|^2)p_j=0 \text{ for all } k+1\leq j\leq n.$$ So $p_j=0$ for all $k+1\leq j\leq n$.
Hence rank $P\leq k$. \qed

\bigskip

\textit{Proof of Theorem \ref{thm2}} \quad 
We first show that it is sufficient to prove the result when $\A$ is a 
subalgebra of block diagonal matrices in $\mat$. If $\A$ is any subalgebra of 
$\mat$ then there exist $n_1, n_2,\ldots,n_k$ with $\sum_i n_i=n$ such that $\A$ 
is $*$-isomorphic to $\oplus_i \mathbb M_{n_i} (\C)$, the $*$-isomorphism 
$\varphi: \A\rightarrow \oplus_i \mathbb M_{n_i} (\C)$ being $\varphi(X)=V^*X V$ 
for some unitary matrix $V\in \mat$ (see \cite[p. 249]{conway}, \cite[p. 74]{davidson}). 
By definition $$\dist(A,\A)=\min_{W\in \A}\|A-W\|.$$ Let $\tilde{A}$ denote the matrix $V^*A V$. 
Since $\|\cdot\|$ is unitarily invariant, we get \begin{equation}
\dist(A,\A)=\dist(\tilde{A},\oplus_i \mathbb M_{n_i} (\C)).\label{21}
\end{equation} 

Next we show that for any density matrix $P$,
\begin{eqnarray}
&&\max\left\{\tr\left(A^*AP-\mathcal C_{\A}(AP)^*\ \mathcal C_{\A}(AP)\ \mathcal C_{\A}(P)^{-1}\right): P\geq 0,\tr P=1 \right\} \nonumber\\
&&\hspace{-1cm}=\max\left\{\tr\left(\tilde{A}^*\tilde{A} \tilde{P}-\mathcal C(\tilde{A} \tilde{P})^*\ \mathcal C(\tilde{A}\tilde{P})\ \mathcal C(\tilde{P})^{-1}\right): \tilde{P}\geq 0,\tr \tilde{P}=1 \right\}, \label{22}
\end{eqnarray}
where $\mathcal C$ is the pinching map as defined in \eqref{pinching}.
Since \begin{equation}
\tr (XY)=\tr (YX),\label{cyclicity of trace}
\end{equation} we have
$$\tr\left(A^*AP-\mathcal C(AP)^*\ \mathcal C(AP)\ \mathcal C(P)^{-1}\right)=\tr\left(V^*A^*APV- V^*\ \mathcal C(AP)^*\ \mathcal C(AP)\ \mathcal C(P)^{-1}\ V\right).$$
Now note that for any $X\in \mat,$ $V^* \mathcal C(X) V=\mathcal C (V^* X V)$. . 
Therefore the above expression is same as 
\begin{equation}
\tr\left(\tilde{A}^*\tilde{A} \tilde{P}-\mathcal C(\tilde{A} \tilde{P})^*\ \mathcal C(\tilde{A}\tilde{P})\ \mathcal C(\tilde{P})^{-1}\right).\label{20a}\end{equation}
This gives \eqref{22}.
So it is enough to prove \eqref{dist to any subalgebra} when $\A$ is a subalgebra of block diagonal matrices. 
We first show that 
\begin{equation}
\max\left\{\tr\left(A^*AP-\mathcal C(AP)^*\ \mathcal C(AP)\ \mathcal C(P)^{-1}\right): P\geq 0,\tr P=1 \right\}\leq \dist(A,\A)^2.\label{bleq}
\end{equation}
Let $P$ be any density matrix.
Then
$\tr (A^*AP)\leq \|A\|^2.$
Therefore 
\begin{equation}
\tr\left(A^*AP-\mathcal C(AP)^*\ \mathcal C(AP)\ \mathcal C(P)^{-1}\right)\leq \|A\|^2.\label{20}
\end{equation}
Let $B\in \A$.  Applying the translation $A\rightarrow A+B$ in \eqref{20} we get 
\begin{equation}
\tr\left((A+B)^*(A+B) P-\mathcal C((A+B) P)^*\ \mathcal C((A+B) P)\ \mathcal C(P)^{-1}\right)\leq \|(A+B)\|^2.\label{AtoA+B}
\end{equation}
We show that the expression on the left hand side  is invariant under this translation.
By expanding the expression on the left hand side of \eqref{AtoA+B}, we get
\begin{eqnarray}
\left(\tr\left(A^*AP-\mathcal C(AP)^*\ \mathcal C(AP)\ \mathcal C(P)^{-1}\right)\right) + \left(\tr\left(B^*AP-\mathcal C(BP)^*\ \mathcal C(AP)\ \mathcal C(P)^{-1}\right)\right)&& \nonumber\\
\hspace{-1cm}+\left(\tr\left(A^*BP-\mathcal C(AP)^*\ \mathcal C(BP)\ \mathcal C(P)^{-1}\right)\right)+\left(\tr\left(B^*BP-\mathcal C(BP)^*\ \mathcal C(BP)\ \mathcal C(P)^{-1}\right)\right).&&\label{bexpanding}
\end{eqnarray}
We show that except for the first term, $\left(\tr\left(A^*AP-\mathcal C(AP)^*\ \mathcal C(AP)\ \mathcal C(P)^{-1}\right)\right)$, the rest of the terms in \eqref{bexpanding} are zero. We shall prove that the second term 
\begin{equation}
\tr\left(B^*AP-\mathcal C(BP)^*\ \mathcal C(AP)\ \mathcal C(P)^{-1}\right).\label{bsecondterm}
\end{equation}
 in \eqref{bexpanding} is zero. The proof for the other two terms  is similar.

By using \eqref{bimodule}, the expression in \eqref{bsecondterm} is equal to 
\begin{equation*}
\tr B^*\left(\mathcal C(AP)- \mathcal C(P)\ \mathcal C(AP)\ \mathcal C(P)^{-1}\right).
\end{equation*}
By \eqref{cyclicity of trace} this is equal to
\begin{equation}
\tr B^*\ \mathcal C(AP)\left(I- \mathcal C(P)^{-1} \mathcal C(P)\right).\label{secondterm}
\end{equation}
If $\mathcal C(P)$ is invertible then this is clearly zero. 
So let $\mathcal C(P)$ be not invertible. This means that if $\mathcal C(P)=\left[\begin{array}{ccc} P_1&\ &\ \\  \ &\ddots&\ \\ \ &\ &P_k\end{array}\right]$, then there exists $i, 1\leq i\leq k,$ such that $P_i$ is not invertible.
Let $U$ denote the block diagonal unitary matrix 
\begin{equation}
U=\left[\begin{array}{ccc} U_1& \ & \ \\
\ & \ddots &\ \\
\ &\ &U_k\end{array}\right],\label{defn of U}
\end{equation} where $U_i=I,$ if $P_{i}$ is invertible and $U_i^* P_{i} U_i=\left[\begin{array}{ccc} \Lambda_i& \ \\
\ & O \\
\end{array}\right],$ if $P_{i}$ is not invertible. (Here $\Lambda_i$ is the diagonal matrix with eigenvalues of $P_{i}$ as its diagonal entries.)
Let $X'$ denote the matrix $U^* X U$. 
Then from \eqref{bimodule} and \eqref{cyclicity of trace}, we get that the expression in \eqref{secondterm} is same as 
\begin{equation}
\tr B'^*\ \mathcal C(A'P')\left(I-  \mathcal C(P')^{-1}\ \mathcal C(P')\right).\label{32}
\end{equation}
Now $\mathcal C(P')=\left[\begin{array}{ccccc} \Lambda_1& \ & \ &\ &\ \\
\ & O &\ &\ &\ \\
\ & \ &\Lambda_2 &\ &\ \\
\ & \ &\ & O &\ \\
\ & \ &\ &\ & \ddots
\end{array}\right]$. Write $A'$ and $P'$ as $2k$-block matrices, $A'=(A'_{rs})_{r,s=1,\ldots,2k} \text{ and } P'=(P'_{rs})_{r,s=1,\ldots,2k},$ respectively such that whenever $P_i$ is not invertible, we have $P_{2i-1, 2i-1}'=\Lambda_i$ and $P
_{2i, 2i}'=O$.

The $(r,r)$-entry of $A'P'$ is $\displaystyle\sum_{s=1}^{2k} A'_{rs} P'_{sr}$. Suppose $P'_{rr}=O$.
Since $P'\geq 0$, we have $P'_{rs}=P'_{sr}=O$ for all $s=1,\ldots,2k$. Hence the $(r,r)$-entry of $A'P'$ is zero.
So let $P'_{rr}\neq O$. Then the $(r,r)$-entry of $\left(I-\mathcal C(P')^{-1} \mathcal C(P')\right)$ is zero.
Therefore the expression in \eqref{32} is zero, and hence the expression in \eqref{secondterm} is zero. 
Therefore from \eqref{AtoA+B}, we obtain
$$\tr\left(A^*A P-\mathcal C(A P)^*\ \mathcal C(A P)\ \mathcal C(P)^{-1}\right)\leq \|(A+B)\|^2,$$
for all $B\in \A$ and for all density matrices $P$.
Equation \eqref{bleq} now follows from here.

To show equality in \eqref{bleq}, let $\dist(A, \A)=\|A_0\|, $ where $A_0=A-B_0$ for some $B_0\in \A$. Then $A_0$ is orthogonal to $ \A$. By Theorem \ref{thm1} there exists a density matrix $P$ such that \begin{equation}
A_0^*A_0 P=\|A_0\|^2 P \label{consequence1b}
\end{equation} and \begin{equation}
\mathcal C(A_0 P)=0, \text{ that is, } \mathcal C(AP)=\mathcal C(B_0 P). \label{consequence2b}
\end{equation}

From \eqref{consequence1b} we get that
 \begin{eqnarray*}
\|A_0\|^2&=& \tr (A-B_0)^*(A-B_0) P \nonumber\\
&=& \tr (A^*AP)-\tr (B_0^* AP)-\tr (A^* B_0 P)+ \tr(B_0^* B_0 P).
\end{eqnarray*}
By using \eqref{bimodule}, we obtain 
\begin{equation}
 \|A_0\|^2=\tr (A^*AP)-\tr (B_0^* \ \mathcal C(AP))-\tr ( B_0 \ \mathcal C(AP)^*)+ \tr(B_0^*\ \mathcal C(B_0P)). \label{31}
\end{equation}

Substituting \eqref{consequence2b} in \eqref{31} we get
\begin{equation}
\|A_0\|^2=\tr(A^*AP)-\tr(B_0^*B_0 P).\label{35}
\end{equation}
Now consider $\tr\left(\mathcal C(AP)^*\ \mathcal C(AP)\ \mathcal C(P)^{-1}\right)$. From \eqref{consequence2b} we see that this is same as $\tr \left(B_0^* B_0\ \mathcal C(P) \mathcal C(P)^{-1} \mathcal C(P)\right)$. If $\mathcal C(P)$ is invertible, then this is equal to $\tr (B_0^* B_0 P)$. If $\mathcal C(P)$ is not invertible, then we define $U$ as done in \eqref{defn of U}. From \eqref{bimodule} and \eqref{cyclicity of trace}, we obtain
 $$\tr \left(B_0^*\ B_0\ \mathcal C(P)\ \mathcal C(P)^{-1}\ \mathcal C(P)\right)=\tr \left(B_0'^*\ B_0'\ \mathcal C(P')\ \mathcal C(P')^{-1}\ \mathcal C(P')\right).$$
 By definition of $U$, this is equal to $\tr \left(B_0'^*\ B'_0\ \mathcal C(P')\right)$, which again by \eqref{bimodule} and \eqref{cyclicity of trace}, is same as $\tr \left(B_0^*\ B_0\ \mathcal C(P)\right)$.
Therefore from \eqref{35} we have
$$\dist(A, \A)^2=\|A_0\|^2= \tr\left(A^*A P-\mathcal C(A P)^*\ \mathcal C(A P)\ \mathcal C(P)^{-1}\right).$$

\section{Remarks}
\begin{enumerate}

\item It is clear from the proof of Theorem \ref{thm1} that  the condition $A^*AP=\|A\|^2 P$ can be replaced by the weaker condition $\tr(A^*AP)=\|A\|^2 $ in the statement of Theorem \ref{thm1}.

\item As one would expect,  the set $\{A: \|A+W\|\geq \|A\| \text{ for all }W\in \W\}$ need not be a subspace.
As an example consider the subspace $\W=\C I$ of $\mathbb M_3(\C)$. Let $A_1=\left[\begin{array}{ccc}
0 &1&0 \\
1 & 0 &1 \\
0 &1&0
\end{array}\right]$ and $A_2=\left[\begin{array}{ccc}
0&0 &1 \\
0 & 0 &0 \\
1 &0 &0
\end{array}\right]$. It can be checked from Theorem \ref{thm1} that $A_1, A_2$ are orthogonal to $\W$. (Take $P=\left[\begin{array}{ccc}
0 &0 &0 \\
0 & 1 &0\\
0 &0 &0
\end{array}\right]$ for $A_1$ and $P=\left[\begin{array}{ccc}
1 &0 &0 \\
0 & 0 &0 \\
0 &0 &0
\end{array}\right]$ for $A_2$, respectively.)
Then $A_1+A_2=\left[\begin{array}{ccc}
0 &1&1 \\
1 & 0 &1 \\
1 &1&0
\end{array}\right]$, and $\|A_1+A_2\|=2$. But $\left\|A_1+A_2-\frac{1}{2}I\right\|=\frac{3}{2}<\|A_1+A_2\|$.
Hence $A_1+A_2$ is not orthogonal to $\W$. 

\item Let $\W=\{X:\tr X=0\}$. Then $ \W^{\perp}=\C I$. In Section 1, we stated that if $A\in \W^{\perp}$ such that $\tr(A^*A)=\|A\|^2$ then $A$ is orthogonal to $\W$. Therefore all the scalar matrices are orthogonal to $\W$.  We show that  if $A\notin \C I$ then there exists a matrix $W$ with $\tr W=0$ such that $\|A+W\|<\|A\|$. Let $\mathcal D A$ and  $\mathcal O A$ denote the diagonal and  off-diagonal parts of $A$, respectively. Then $\mathcal O A\in \W$, $A-\mathcal O A=\mathcal D A$ and $\|\mathcal D A\|\leq \|A\|$. So it is enough to find $W\in \W$ such that $\|\mathcal D A+W\|<\|\mathcal D A\|$.
Let $\mathcal D A=\diag\left(a_1,\ldots,a_1,a_2,\ldots,a_2,\ldots,a_k,\ldots,a_k\right),$
where each $a_j$ occurs on the diagonal $n_j$ times and $n_1+\cdots+n_k=n$. Assume $\|\mathcal D A\|=1$. 
Take $W=\diag\left(\frac{a_2-a_1}{k n_1},\ldots,\frac{a_2-a_1}{k n_1},\frac{a_3-a_2}{k n_2},\ldots,\frac{a_3-a_2}{k n_2},\ldots,\frac{a_{k}-a_{k-1}}{k n_{k-1}},\ldots,\frac{a_{k}-a_{k-1}}{k n_{k-1}},\right.\\ \left.\frac{a_1-a_k}{k n_k},\ldots, \frac{a_1-a_k}{k n_k}\right).$
Then $W$ has trace zero and $\mathcal D A+W=\diag\left(\frac{(n_1-1)a_1+a_2}{n_1},\ldots,\right.\\ \left.\frac{(n_1-1)a_1+a_2}{n_1},\right.$ $\left. \frac{(n_2-1)a_2+a_3}{n_2},\ldots,\frac{(n_2-1)a_2+a_3}{n_2},\ldots,\frac{(n_k-1)a_{k-1}+a_k}{n_k},\ldots,\frac{(n_k-1)a_{k-1}+a_k}{n_k}\right).$ It is easy to check that $\|\mathcal D A+W\|<1$. Hence for this particular $\W$ we have that $\{A: \|A+W\|\geq \|A\| \text{ for all } W\in \W\}=\W^{\perp}=\C I$.


\item In Theorem \ref{thm1}, $m(A)$ is the  best possible upper bound on rank $P$. Consider $\W=\{X:\tr X=0\}$. From Remark 2, we get that if a matrix $A$ is orthogonal to $\W$ then it has to be of the form $A=\lambda I$, for some $\lambda \in \C$. When $A\neq 0$ then $m(A)=n$. Let $P$ be any density matrix satisfying $AP\in \W^{\perp}$. Then $AP=\mu I$, for some $\mu \in \C, \mu\neq 0$. If $P$ also satisfies  $A^*AP=\|A\|^2 P$, then we get $P=\frac{\mu}{\lambda} I$. Hence rank $P=n=m(A)$.

\item For $n=2$ and $\A$ any subalgebra of $\mathbb M_2(\C)$, we can restrict maximum on the right hand side of \eqref{dist to any subalgebra} over rank one density matrices.  By the same argument as in the proof of Theorem \ref{thm2} it is sufficient to prove this for $\D_2(\C)$, the subalgebra of diagonal matrices with complex entries. We show
\begin{equation}
\dist(A,\D_2(\C))^2=\max_{\|x\|=1} \left(\|Ax\|^2-\tr \Delta(Axx^*)^*  \Delta(Axx^* )\Delta(xx^*)^{-1}\right),
\end{equation}
where $\Delta$ is the projection onto $\D_2(\C)$. From Theorem \ref{thm2} we have
$$\max_{\|x\|=1} \left(\|Ax\|^2-\tr \Delta(Axx^*)^*  \Delta(Axx^* )\Delta (xx^*)^{-1} \right)\leq \dist(A,\D_2(\C))^2.$$
Note that
\begin{equation}
\dist(A, \D_2(\C))\leq \|\mathcal O A\|.\label{38}
\end{equation}
Let $A=\left[\begin{array}{ccc}
a& b\\
c & d
\end{array}\right]$
and without loss of generality assume that $|b|\geq |c|$. Then $\|\mathcal O A\|=|b|$. For $x=\left[\begin{array}{ccc} 0\\ 1 \end{array}\right]$
$$\|Ax\|^2-\tr \Delta(Axx^*)^*  \Delta(Axx^* )\Delta (xx^*)^{-1} =\|\mathcal O A\|^2.$$ 
Combining this with \eqref{38}, we obtain
$$\dist(A,\D_2(\C))^2\leq \max_{\|x\|=1} \left(\|Ax\|^2-\tr \Delta(Axx^*)^*  \Delta(Axx^* )\Delta (xx^*)^{-1}\right).$$

\bigskip

 \item For $n=2$ and $\A$ any subalgebra of $\mathbb M_2(\C)$, we note that $$\{A: \|A+W\|\geq \|A\| \text{ for all } W\in \A\}=\A^{\perp}.$$ Again it is enough to show that
$$\{A: \|A+W\|\geq \|A\| \text{ for all } W\in \D_2(\C)\}=\D_2(\C)^{\perp}.$$
If $A$ is an off-diagonal $2\times 2$ matrix, that is, $A=\left[\begin{array}{ccc}
0 & b\\
c & 0 
\end{array}\right]$ then by Theorem 2.1 in \cite{bhatiachoidavis} we obtain $\|A+W\|\geq \|A\| \text{ for all } W\in \D_2(\C)$.
Conversely let $A\in \mat$ be such that $\|A+W\|\geq \|A\| \text{ for all } W\in \D_2(\C)$. Then by taking $W=-\mathcal D A$, we have $A+W=\mathcal O(A)$.  Again by using Theorem 2.1 in \cite{bhatiachoidavis} we obtain that $\|\mathcal O(A)\|= \|A\|$. So $A$ is of the form $\left[\begin{array}{ccc}
a & b\\
c & d 
\end{array}\right]$, where $\|A\|=\max\{|b|,|c|\}.$ Since norm of each row and each colum is less than or equal to $\|A\|$, we get that $a=d=0$. Hence $A\in \D_2(\C)^{\perp}$. 
\end{enumerate}
\bigskip

\textbf{Acknowledgement.} {I would like to thank Professor Rajendra Bhatia for several useful discussions and Professor Ajit Iqbal Singh for helpful comments in this paper. }

\end{document}